\documentclass[11pt]{article}

\def\R{\ifmmode{\rm I\mkern-3.1mu
R\mkern1mu}\else{\rm I\kern-.18em
R\hskip1pt\ }\fi\relax}

\def\Z{\ifmmode{\it Z\mkern-7.8mu
Z\mkern2mu}\else{\it Z\kern-.28em
Z\hskip1pt\ }\fi\relax}

\def\Q{\ifmmode{\rm Q\mkern-10mu
l\mkern4.5mu}\else{\rm Q\kern-.57em
l\hskip3pt\ }\fi\relax}

\def\N{\ifmmode{\rm I\mkern-3.1mu
N\mkern0.5mu}\else{\rm I\kern-.16em
N\hskip0.5pt\ }\fi\relax}

\def\C{\ifmmode{\rm C\mkern-8.8mu
l\mkern4mu}\else{\rm C\kern-.48em
l\hskip2.6pt\ }\fi\relax}

\def\Hyper{\ifmmode{\rm I\mkern-3.1mu
H\mkern0.5mu}\else{\rm I\kern-.16em
H\hskip0.5pt\ }\fi\relax}

\title{A theorem of Poincar\'e-Hopf type}
\author{St\'ephane Simon}
\date{\today}

\newtheorem{defi}{Definition}[section]

\newtheorem{lem}{Lemma}[section]
\newtheorem{theo}{Theorem}[section]
\newtheorem{corol}{Corollary}[section]
\newtheorem{prop}{Proposition}[section]

\begin{document}

\bibliographystyle{alpha}
\maketitle

\begin{abstract}
We compute (algebraically) the Euler
characteristic of a complex of sheaves with constructible
cohomology. A stratified Poincar\'e-Hopf formula is then a
consequence of the smooth Poincar\'e-Hopf theorem and of
additivity of the Euler-Poincar\'e characteristic with compact
supports, once we have a suitable definition of index.
\end{abstract}
{\small
}

\section{Introduction}

M.-H. Schwartz has defined radial vector fields in \cite{Sc1} and
extended the classical Poincar\'e-Hopf theorem to real analytic
sets, equip\-ped with a Whitney stratification for these vector
fields \cite{Sc4}, \cite{Sc5}. In their turn, H. King and D.
Trotman have extended M.-H. Schwartz's result to more general
singular spaces and generic vector fields \cite{KT}.

Radial \cite{Sc1}, \cite{Sc3} (and totally radial, see \cite{KT},
\cite{Si2}) vector fields are important because of their relation
with Chern-Schwartz-Mac Pherson classes. Chern-Mac Pherson classes
are written as an integral combination of Mather classes of
algebraic varieties with coefficients determined by local Euler
obstructions \cite{Mac1}. A transcendental definition (and the
original one) of local Euler obstruction is the obstruction to
extend a lift of a radial vector field, prescribed on the link of
a point in the base, inside a whole neighborood of Nash transform.
Chern-Schwartz classes \cite{Sc1}, \cite{Sc3} (which lie in
cohomology of the complex analytic variety) are defined as the
obstruction to extend a radial frame field given on a sub-skeleton
of a fixed triangulation. These two points of view coincide:
Chern-Mac Pherson classes are identified with Chern-Schwartz
classes by Alexander duality \cite{BS}. In \cite{BBFGK}, it is
shown that these Chern-MacPherson-Schwartz classes can be realised
(in general not uniquely) in intersection homology with middle
perversity.

This paper concerns a Poincar\'e-Hopf theorem in intersection
homology for a stratified pseudo-manifold $A$ (\cite{GM1}) and a
vector field $v$ which does not necessarily admit a globally
continuous flow. Our main result is that we still have a
Poincar\'e-Hopf formula when the vector field is semi-radial
\cite{KT}~: $$I\chi_c^{\overline p}(A)=\sum_{v(x)=0}Ind^{\overline
p}(v,x).$$ More precisely, we compute (algebraically) the Euler
characteristic of a complex of sheaves with constructible
cohomology. A stratified Poincar\'e-Hopf formula is then a
consequence of the smooth Poincar\'e-Hopf theorem and of
additivity of the Euler-Poincar\'e characteristic with compact
supports, once we have a suitable definition of index.\\ Given a
vector field with isolated singularities on a singular space,
which admits a globally continuous flow, one can already deduce a
Poincar\'e-Hopf theorem from a Lefschetz formula in intersection
homology with middle perversity \cite{GM2}, \cite{GM3},
\cite{Mac2}.

A. Dubson announced in \cite{Du} a  formula similar
to ours for a constructible complex in a complex analytic
framework. In \cite{BDK}, J.-L. Brylinski, A. Dubson and M.
Kashiwara expressed the ``local characteristic'' of a holonomic
module as a function of multiplicities of polar varieties and
local Euler obstructions.

M. Goresky and
R. MacPherson have proved a Lefschetz fixed point theorem for a subanalytic morphism
and
constructible
complex of sheaves \cite{GM3}. They show that a weakly hyperbolic
morphism
({\it i.e.} whose fixed points are {\sl weakly hyperbolic}) can be lifted
to a morphism (not necessarily unique) at the level of sheaves.
The Lefschetz number can be written
as a sum of contributions of the various connected components of fixed points, a
component being itself possibly stratified; every contribution is a sum of
multiplicities (relative to the morphism), weighted by Euler
characteristics in compactly supported cohomology of the strata of the connected
component.
%For a flow, we retrieve of course the Euler characteristic of the
%stratified set and our result in identifying the Lefschetz index with
%vector index.

% Another result on the Euler characteristic is given by S. Cappell and J. Shaneson who
% give in \cite{CS} formulae
% (inductive on strata, and verifying a ``stratified multiplicative property'') for
% the Euler characteristic of a complex projective variety (possibly singular),
% under the
% effect of a morphism.
% More precisely, the characteristic of the source is a sum of the
% product of the target characteristic by that of the generic fibre, with
% coefficients involving the fibre of the projective completion.

In Section 2 we give a formula to calculate the characteristic of
a constructible complex of sheaves. Then, in section 3, we apply
the preceding results to the intersection chain complex. A brief
recall of definitions and results on stratified vector fields is
given in section 4. A theorem of Poincar\'e-Hopf type appears in
section 5, where the vector field considered is totally (or only
semi-) radial. Sections 6 and 7 are devoted to illustrate the
theorems of section 4.

I am very grateful to D. Trotman for numerous valuable conversations, and for
his constant encouragement. I thank equally J.-P.
Brasselet for the discussions we have had related to this
work. I am also greatly indebted to E. Leichtnam for his active
interest and for many helpful suggestions during the preparation of
the paper.

\section{A formula to calculate the Euler-Poincar\'e characteristic of a complex
of sheaves with constructible cohomology}

First we recall some definitions. Let $R$ be a
principal ideal domain. We shall consider sheaves of $R-$modules.

\begin{defi} A {\rm stratified set} $A$ is
a topological space which is a union of a locally finite family of
disjoint, connected subsets (strata) which are smooth manifolds,
satisfying the frontier condition. We shall denote by $\cal A$ the
set of strata and suppose that this stratification is fixed once
and for all.
\end{defi}
\begin{defi} Let $A$ be a stratified set. We say that
$A$ is {\rm compactifiable} if there exists a compact
abstract stratified set $(B,\cal B)$ {\rm(\cite{Mat1},
\cite{Mat2}, \cite{T}, \cite{V2})}, such that $A\subseteq B$ is a
locally closed subset of $B$ which is a union of elements of $\cal
B$. We then say that $(B,\cal B)$ is a {\rm compactification} of
$A$.
\end{defi}
\begin{defi} Let $A$ be a stratified set and $\cal F$ a
sheaf on $A$. We say that $\cal F$ is {\rm $\cal A-$constructible} on
$A$ if for every stratum $X$ of $\cal A$, the sheaf
${\cal F}_{\vert X}$ is locally constant of finite rank on $R$.
\end{defi}

Recall that
$H_c(A;{\cal F})\cong\Hyper_c(A;{\cal F})$
where $\Hyper_c$ denotes hypercohomology with compact supports. 
% and value in $\cal F$).
%where we have identified $\cal F$ with the complex of sheaves concentrated in
%degree $0$ and with null differential
As usual, suppose that $H_c^p(A;{\cal F})$ has finite rank for
$p\geq 0$ and is null for large enough $p$. Then we call {\sl
Euler characteristic of $A$ with compact supports and coefficients
in $\cal F$}, the alternating sum of the ranks of the modules
$H_c^p(A;{\cal F})$ and denote it by $\chi^c(A;{\cal F})$. When
the sheaf $\cal F$ is the constant sheaf $R$, we simply write
$\chi^c(A)$. We shall see that the Euler characteristic is always
defined in our situation.
\begin{prop} Let $X$ be a locally compact topological space,
$\cal G$ a locally constant sheaf on $X$ of finite rank $g$ and
suppose that $X$ admits a finite partition $\cal T$ into open
simplexes, i.e. there exists a finite simplicial complex (resp.
subcomplex, possibly empty) $K$ (resp. $L$) and a homeomorphism
$\varphi:K\setminus L\rightarrow X$. Then $$\chi^c(X;{\cal
G})=\chi^c(X).g.$$
\end{prop}
{\it Proof}. As simplexes of $X$ are contractible, the restriction
of $\cal G$ is isomorphic to the constant sheaf over any one of
them. Consider the finite union $U$ of open simplexes of
ma\-xi\-mal dimension  $m$. By induction on $m$ and using the long
cohomological exact sequence (with compact supports) of $(U,X)$,
we are reduced to showing the result for $U$. But, applying
Mayer-Vietoris to the partition of $U$, this shows that
$\chi^c(U;{\cal G}_{\vert U})$ is well defined and establishes the
formula.
%We are going to do induction on $m$.
%\smallskip
%\noindent
%If $m=0$ then $X$ is a finite union of points and we have
%$\chi^c(X;{\cal G})=\chi^c(X).g$.
%\smallskip
%\noindent
%Suppose now $m>0$. We have a long exact sequence in cohomology~:
%$$\cdots\rightarrow H^p_c(U;{\cal G}_{\vert U})\rightarrow H^p_c(X;{\cal G})
%\rightarrow
%H^p_c(X\setminus U;{\cal G}_{\vert X\setminus U})\rightarrow\cdots$$
%\noindent
%As $U$ is a disjoint union of $N_m$ simplexes $\sigma$ and
%$\cal G$ is constant on all $\sigma$, we see that $H^p_c(U;{\cal G}_{\vert U})$
%is finite for all $p\geq 0$. On the other hand $H^p_c(X\setminus U;
%{\cal G}_{\vert X\setminus U})$ is
%finite by induction. Finally $H^p_c(X;{\cal G})$ is finite for all $p\geq 0$.
%\par
%Applying Mayer-Vietoris to the partition of $U$, we find~:
%$$\begin{array}{ll}
%\chi^c(U;{\cal G}_{\vert U})&=\sum_{k=1}^{N_m}\chi^c(\sigma_k;{\cal
%G}_{\vert\sigma_k})\\
%\null                       &=\sum_{k=1}^{N_m} (-1)^m g\\
%\null &=\chi^c(U).g.
%\end{array}$$
%\noindent
%We have thus :
%$$\begin{array}{ll}
%\chi^c(X;{\cal G})&=\chi^c(X\setminus U;{\cal G}_{\vert X\setminus
%U})+\chi^c(U;{\cal G}_{\vert U})\\
%\null&=(\chi^c(X\setminus U)+\chi^c(U)).g\\
%\null&=\chi^c(X).g.
%\end{array}$$
\begin{prop} Let $A$ be a compactifiable stratified set
and $(B,\cal B)$ a compactification of $A$. Let
${\cal A}=(X_i)_{i\in\{1,\ldots,N\}}$ be the strata of $A$ and $\cal F$
an $\cal A-$constructible sheaf. Then we have~:
$$\chi^c(A;{\cal F})=\sum_{i=1}^N\chi^c(X_i).rk\ {\cal F}_{\vert X_i}.$$
\end{prop}
{\it Proof}. Write $\overline{A}$ for the closure of $A$ in $B$.
Thanks to the triangulation theorem for abstract stratified sets of M. Goresky
\cite{Gor}, there exists a triangulation $\cal T$ of
$\overline{A}$ adapted to the stratification $\overline{\cal A}$. As
$\overline{A}$ is compact, this triangulation is finite. Moreover, it is also
adapted to ${\cal A}$.

\noindent
We are going to do induction on the number of strata of $A$ and
apply  the method of proof of proposition 2.1. Let
$X$ be a stratum of maximal depth (\cite{V2}) in  $A$. Remark that $X$ is closed
in $A$.
\smallskip
\noindent
If $A=X$ we apply proposition 2.1 with $X$ and $\cal F$.

\noindent
Suppose the cardinal of $\cal A$ is strictly greater than $1$.

\noindent
We have then a long exact sequence in cohomology :
$$\cdots\rightarrow H_c^p(A\setminus X;{\cal F}_{\vert A\setminus X})\rightarrow
H_c^p(A;{\cal F})
\rightarrow H_c^p(X;{\cal F}_{\vert X})\rightarrow\cdots$$
As the number of strata of $A\setminus X$ is strictly smaller than
that in $A$, we can apply the induction hypothesis to $A\setminus X$
and ${\cal F}_{\vert A\setminus X}$. This shows that
$rk\ H_c^p(A;{\cal F})$ is finite, so $\chi^c(A;{\cal F})$ is defined.
On the other hand, we have~:
$$\chi^c(A;{\cal F})=\chi^c(A\setminus X;{\cal F}_{\vert A\setminus X})+
\chi^c(X;{\cal F}_{\vert X}).$$
\noindent
We conclude by using the induction hypothesis and proposition 2.1.
\bigbreak
Let ${\cal F}^{\bullet}$ be a complex of sheaves. Let
$\cal H^{\bullet}(F^{\bullet})$  be the complex
of derived sheaves.
\bigbreak
\begin{defi}Let $A$ be a compactifiable stratified set
and ${\cal F}^{\bullet}$ a complex of sheaves on $A$. We say that
${\cal F}^{\bullet}$ has
{\rm $\cal A-$constructible cohomology} if~:
{\parindent=1cm
\item {(i)} ${\cal F}^{\bullet}$ is bounded
\item {(ii)} ${\cal H}^{\bullet}({\cal F}^{\bullet})$
is $\cal A-$constructible.
}
\end{defi}
\begin{theo} Let $A$ be a compactifiable stratified set,
${\cal A}=(X_i)_{i\in\{1,\ldots,N\}}$ its stratification and ${\cal F}^{\bullet}$
a complex of $c-$acyclic sheaves with $\cal A-$constructible cohomology.  Then we
have~:
$$\chi^c(A;{\cal F}^{\bullet})=\sum_{q=-N_1}^{N_2}(-1)^qrk\
\Hyper_c^q(A;{\cal F}^{\bullet})=\sum_{i=1}^N\chi^c(X_i)
\chi({({\cal H}^{\bullet}({\cal F}^{\bullet})_{\vert X_i})}_{x_i})$$
\noindent
\nobreak
where ${\cal F}^p=0$ except for
$-N_1\leq p\leq N_2$ and $x_i$ is any point of $X_i$, $1\leq i\leq N$.
\end{theo}
{\it Proof}. As ${\cal F}^p$ is $c-$acyclic for all $p\in\Z$, we have
$H^p(H^q_c(A;{\cal F}^{\bullet}))=0$ for all $p\in\Z$ and $q\geq 1$. So the
second
spectral sequence,  of second term
$`E_2^{pq}=H^p(H^q_c(A;{\cal F}^{\bullet}))$, degenerates.
As ${\cal F}^\bullet$ is bounded, the filtration of the associated double
complex
is regular, so the first spectral sequence is convergent and we have
according to theorem 4.6.1 of \cite{God} p. 178~:
$$E^{p,q}_2=H^p_c(A;{\cal H}^q({\cal F}^{\bullet}))\Rightarrow
\Hyper^{p+q}_c(A;\cal F^{\bullet}).$$
As ${\cal F}^{\bullet}$ has $\cal A-$constructible cohomology  and $A$ is
compactifiable, we can define~:
$$\begin{array}{ll}
\chi(E_2)&=\sum_{p\in\N,q\in\Z}(-1)^{p+q}rg\ H^p_c(A;{\cal H}^q(
{\cal F}^{\bullet}))\\
&=\sum_{q\in\Z}(-1)^q\sum_{p\in\N}(-1)^prg\ H^p_c(A;{\cal H}^q
(\cal F^{\bullet}))\\
&=\sum_{q=-N_1}^{N_2}(-1)^q\chi^c(A;{\cal H}^q(\cal F^{\bullet}))
\end{array}$$
\noindent
for $\cal F^{\bullet}$ is bounded. Remark that, since $A$ is triangulable,
every point of $A$ (which is paracompact) admits a neighborhood
homeomorphic to a subspace of some $\R^p$, so that $A$ is
of cohomological dimension lower or equal to $p$ ($<\infty$ because $A$ is
compactifiable), according to theorem 5.13.1 of \cite{God} p. 237.
\smallbreak
Apply then proposition 2.2 to  $A$ and ${\cal H}^q(\cal F^{\bullet})$~:
$$\begin{array}{ll}\chi(E_2)&=\sum_{q=-N_1}^{N_2}(-1)^q\sum_{i=1}^N\chi^c(X_i)
rg\ ({\cal H}^q({\cal F}^{\bullet})_{\vert X_i})_{x_i}\\
\null&=\sum_{i=1}^N\chi^c(X_i)\sum_{q=-N_1}^{N_2}(-1)^q
rg\ ({\cal H}^q({\cal F}^{\bullet})_{\vert X_i})_{x_i}\\
\null&=\sum_{i=1}^N\chi^c(X_i)\sum_{q=-N_1}^{N_2}(-1)^q
rg\ {\cal H}^q({\cal F^{\bullet}})_{x_i}\\
\null&=\sum_{i=1}^N\chi^c(X_i)\chi({\cal H}^{\bullet}({\cal F}^{\bullet})
_{x_i})
\end{array}$$
\noindent
with $x_i\in X_i$ for $i\in\{1,\ldots,N\}$.
\smallskip
\noindent
As $E_{r+1}=H(E_r)$, we have $\chi(E_{r+1})=\chi(E_r)$ for all $r\geq 2$. So
$\chi(E_r)=\chi(E_2)$ for all $r\geq 2$.

\noindent
As $\cal F^{\bullet}$ is bounded, $E^{p,q}_2=0$ for $q$ big enough or small
enough and $p\in\N$. Thus the spectral sequence degenerates and so
$$E^{p,q}_r=E^{p,q}_{\infty}$$
\noindent
\nobreak
for $r$ big enough.
\goodbreak

\noindent
Hence
$$\chi(E_{\infty})=\chi(E_r)=\chi(E_2).$$
But $(E^{p,q}_{\infty})_{p+q=s}$ is the associated graded module to
$\Hyper^s_c(A;{\cal F}^{\bullet})$. We have thus~:
$$rg\ \Hyper^s_c(A;{\cal F}^{\bullet})
=\allowbreak\displaystyle\sum_{p+q=s}rg\ E^{p,q}_{\infty}.$$
\goodbreak
\noindent
Finally
$$\begin{array}{ll}\chi^c(A;{\cal F}^{\bullet})&=\sum_{s\in\Z}(-1)^srg\
\Hyper^s_c(A;{\cal F}^{\bullet})\cr
\null&=\chi(E_{\infty})\cr
\null&=\chi(E_2)\cr
\null&=\sum_{i=1}^N\chi^c(X_i)\chi({\cal H}^{\bullet}({\cal F}
^{\bullet})_{x_i}).
\end{array}$$
\bigbreak \noindent {\it Remark}. Theorem 2.1 works also with the
weaker hypothesis of (finite) triangulabity.

\section{Application to intersection homology}

Suppose now that $A$ is a pseudo-manifold, and let $\cal A$ be its
stratification. Here the strata of $A$ will no longer be necessarily connected,
but we shall work with connected components of strata. We denote by
$L_{x}$ the link of the point $x$ in $A$.
%If A is also an abstract stratified set then the
%topology of $L_{x}$ is constant on the (unique component of the) stratum which
%contains $x$ (\cite{Mat1}), otherwise this results from the definition
%of a pseudo-manifold.
\begin{prop}[\cite{Bo}] Let $A$ be an
$n$-pseudo-manifold and $\overline p$ a perversity. Let
$IC_{\bullet}^{\overline p}$ be the intersection chain complex for
perversity $\overline p$ with coefficients in $R$ {\rm\cite{GM3}}
and set ${\cal IC}^{\bullet}_{\overline p}=$sheaf associated to
the presheaf $\{U\mapsto IC^{\overline p}_{n-\bullet}(U)\}$. Then
${\cal IC}_{\overline p}^{\bullet}$ is a complex of $c-$soft
sheaves (so $c-$acyclic). Moreover we have~:
$$\Hyper_c^{\bullet}(A;{\cal IC}^{\bullet}_{\overline p})=
IH_{n-\bullet}^{\overline p}(A;R).$$
\end{prop}
\begin{prop}[Proposition 2.4 of \cite{GM1}] Let $A$ be an
$n$-pseudo-manifold, $x$ any
point in a stratum
$X^k$ of $A$ of dimension $k$ and $L_x$ the link of $X^k$ at $x$ in
$A$. The fibre of the
complex of derived sheaves
${\cal H^{\bullet}(IC}_{\overline p}^{\bullet})$ is given by~:
$${\cal H}^i({\cal IC}^{\bullet}_{\overline p})_x=\cases{
\cases{
IH^{\overline p}_{n-i-k-1}(L_x) & if $i\leq p_{n-k}$\cr
0 & otherwise\cr}
& if $x\in X^k\subset A\setminus A_{reg}$\cr
\cases{ R & if $i=0$\cr
0 & otherwise\cr}
& if $x\in X^n\subset A_{reg}$.\cr}
$$
\end{prop}
As usual the {\it Euler-Poincar\'e characteristic
in intersection homology} $I\chi^{\overline p}_c(A)$ of an $n$-pseudo-manifold $A$ is the
Euler-Poincar\'e characteristic with compact supports of the complex of sheaves
${\cal IC}_{\overline p}^{\bullet}$ multiplied  by $(-1)^n$, {\it i.e.}
$I\chi^{\overline p}_c(A)=(-1)^n\chi^c(A;{\cal IC}_{\overline p}^{\bullet}).$
\bigskip
\begin{theo} Let $A$ be an $n$-pseudo-manifold such that
$(A,\cal A)$ is compactifiable,
$N$ the
number of connected components of strata of $A$ and
$\overline p$ a perversity.\\
We have~:
$$I\chi_c^{\overline p}(A)=\sum_{i=1}^N(-1)^n\chi^c(X_i)
\sum_{j=0}^{p_{n-\dim X_i}}
(-1)^j rg\ IH^{\overline p}_{n-j-\dim X_i-1}(L_{x_i};R)$$
\noindent
where $x_{i}$ is an arbitrary point of $X_{i}$ for $1\leq i\leq N$ and
we make the
convention that
$rg\ IH^{\overline p}_{-1}(L_{x_i};R)=1$ if
$\dim X_i=n$.
\end{theo}
{\it Proof}. Application of theorem 2.1 and proposition 3.2.
\bigbreak
\noindent
{\it Remark}. As in theorem 2.1 we can weaken
the hypothesis by only assuming the existence of a (finite) triangulation
compatible with the stratification.
\begin{prop} Let $A$ be a $2n$-pseudo-manifold such that
$(A,\cal A)$ is compactifiable, the dimension of strata being even and let
$\overline m$ be the middle perversity. We have~:
$$I\chi^{\overline m}_c(A)=\sum_{i=0}^{n}\sum_{k=1}^{N_i}\chi_c(X_k^{2i})
\sum_{j=0}^{n-i-1}
(-1)^jrg\ IH^{\overline m}_{2n-j-2i-1}(L_{x^i_k};R)$$
\noindent
where we have written $X_k^{2i}$ (resp. $N_i$) for the $k$-th connected component of
the stratum (resp. number of connected components of the stratum) of
dimension $2i$, $x_{i}^k$ an arbitrary point of $X_k^{2i}$ and
$\displaystyle\chi_c(X)=\sum_{i=0}^{\dim X}(-1)^irg\ H^c_i(X;R)$.
\end{prop}
{\it Proof}. We apply theorem 2.1 with
$\overline p=\overline m$ and we remark that $\chi_c(X)=\chi^c(X)$ for
a manifold $X$ of even dimension.
\section{Totally radial and semi-radial vector fields on abstract stratified
sets}

M.-H. Schwartz constructed certain frame fields to define (by obstruction)
her Chern-Schwartz classes in the cohomology of a singular complex analytic
variety equipped
with a Whitney stratification \cite{Sc1},
\cite{Sc3}. These
were
called radial fields.
When one is concerned with 1-frame fields (i.e. vector fields), they are
called radial vector fields. She showed that they verified a Poincar\'e-Hopf
formula \cite{Sc4}, \cite{Sc5}.

This section is an easy transcription to abstract stratified sets
of some notions and results of \cite{KT} which were given in the
more general setting of  ``mapping cylinder stratified space with
boundary''. In their paper, H. King and D. Trotman extend M.-H.
Schwartz's work on Poincar\'e-Hopf formulas, to more general
spaces, and to generic vector fields. Notice that abstract
stratified sets are not (necessarily) embedded nor are vector
fields (necessarily) continuous.
\begin{defi}[\cite{KT}] Let $(A,\cal A)$ be an abstract stratified set and
$v$ a stratified vector field on $A$ {(\rm
\cite{Mat1},\cite{Mat2}, \cite{T}, \cite{V2}}). We say that $v$ is
a {\rm totally radial} vector field if for all strata
$X\in\cal A$ there exists a neighborhood $U_X$ of $X$ in the
control tube $T_X$ such that $d\rho_X(v)>0$ on $U_X\setminus X$
(i.e. $v$ is pointing outwards with respect to the level
hypersurfaces of the control function $\rho_{X}$).
\end{defi}
In \cite{KT} such a vector field was called {\sl radial}. To avoid confusion with
the radial vector fields of M.-H. Schwartz, we have adopted the terminology
{\sl totally radial}, which also expresses the fact that one imposes that
$d\rho_X(v)>0$ on a whole neighborhood $U_X$ of $X$ in $T_X$. The analoguous
condition is only imposed on a neighbourhood of some closed subset of $X$ by
M.H. Schwartz. See \cite{Si2} for a detailed discussion of the differences
between the radial fields of \cite{Sc4}, \cite{Sc5} and the radial fields
of \cite{KT}, called
totally radial here.
\begin{prop} Let $(A,\cal A)$ be an abstract stratified set and $Y$
a stratum of
$A$. Then there exists a vector field $\xi_Y$ on $T_Y\setminus Y$ such that~:
$$\hbox{for all }y\hbox{ in }T_Y\setminus Y\hbox{ we have }
\cases{{{\rho_Y}_*}_y(\xi_Y)=1\cr
{{\rho_X}_*}_y(\xi_Y)=0&if $X<Y$.\cr}$$
\end{prop}
{\it Proof}. It suffices to consider the stratified submersion
$(\pi_Y,\rho_Y):T_Y\setminus Y\rightarrow Y\times\R_+^*$ and to lift the
constant field $(0,\partial_t)$ to a field $\xi_Y$ on $T_Y\setminus Y$. Thanks
to the compatibility
conditions, we see that ${\rho_X}_*(\xi_Y)=0$ for $X<Y$.
\begin{defi}[\cite{KT}] Let $(A,\cal A)$ be an abstract stratified set,
$v$  a stratified vector field on $A$ and $Y$ a stratum of $A$. Let
$(Y_i)_{1\leq i\leq m}$ be the strata
such that  $Y<Y_i$. Set
$B_Y(v)=\{x\in T_Y\setminus Y\vert(\exists c_i\in\R_-\vert 0\leq i\leq m):
v(y)=c_0\xi_Y(y)+\sum_{j=1}^mc_i\xi_{Y_i}(y)\hbox{ with } c_0<0\}$.
A point
$x\in\overline{B_Y}(v)\cap Y$ is called a {\rm virtual zero } of $v$.
\end{defi}
\begin{defi}[\cite{KT}] Let $(A,\cal A)$ be an abstract stratified set and
$v$ a stratified vector field on $A$. Then $v$ is called {\rm semi-radial} if
$v$ has no virtual zero.
\end{defi}
\noindent
{\it Examples}. Totally radial vector fields, and controlled vector fields,
are semi-radial.\\
% \makebox[10cm][c]{
%   \vbox to 3cm{
%     \vfill
%     \special{pict=test.pict}
%   }
%  %\begin{center}
%  {\sl A totally radial vextor field}
%  %\end{center}
% }

\begin{defi} Let $(A,\cal A)$ be a compactifiable stratified
set, $(B,\cal B)$ a compactification of
$A$ such that $\cal A\subseteq B$ and
$v$ a stratified vector field on $A$. We say that $v$ is {\rm strongly
totally radial} (resp. {\rm strongly semi-radial}) if and only if there exists a
totally radial (resp. semi-radial) extension  $u$ of $v$ to
$(B,\cal B)$.
\end{defi}
\begin{lem}[\cite{KT}] Let $(A,\cal A)$ be an abstract stratified set (resp.
compactifiable stratified set) and $v$ a
semi-radial (resp. strongly semi-radial) vector field with isolated singularities
on $A$. Then there exists
a (resp. strongly) totally radial vector field $v'$ having the same singularities
as
$v$ and the same
indices at these points.
\end{lem}
\section{Towards a Poincar\'e-Hopf theorem}
\begin{defi} Let $A$ be an n-pseudo-manifold,
$\overline p$ a perversity and $x$ a point of a stratum $X$.
We call {\rm multiplicity} of $A$ at $x$ for perversity
$\overline p$ the following integer~:
$$m^{\overline p}_x(A)=\cases{\displaystyle\sum_{i=n-p_{n-\dim X}}^n(-1)^irg\
IH^{\overline p}_{i-\dim X-1}(L_x;R)
& if $x\in A\setminus A_{reg}$\cr
(-1)^n & if $x\in A_{reg}$.\cr}$$
\end{defi}
\noindent
{\it Remark}. The multiplicity is nothing else than
$I\chi_c^{\overline p}(A,A-\{x\})$ (which equals
$(-1)^n$ if $x\in A_{reg}$).
\begin{defi}
Let $A$ be an n-pseudo-manifold such that $(A,\cal A)$ is a
compactifiable abstract stratied set, $\overline p$ a perversity
and $v$ a stratified vector field having an isolated singularity
at $x\in X$. We call {\rm singular index} of $v$ at $x$, and we
denote by $Ind^{\overline p}(v,x)$ the integer: $$Ind^{\overline
p}(v,x)=m^{\overline p}_x(A).Ind(v,x).$$ \noindent Recall that if
the stratum $X$ is reduced to a point, then $Ind(v,x)=1$.
\end{defi}
\begin{theo} Let $A$ be an $n$-pseudo-manifold such that
$(A,\cal A)$ is a compacti\-fia\-ble abstract stratified set, $\overline p$
a perversity and $v$ a strongly semi-radial vector field admitting a finite number
of singularities on $A$. We have~:
$$I\chi_c^{\overline p}(A)=\sum_{v(x)=0}Ind^{\overline p}(v,x).$$
\end{theo}
{\it Proof}. As in \cite{Be}, for all strata $X$ of $A$, let $f_X$
be a carpeting function, i.e. let $U_{b(X)}$ be a neighborhood of
$b(X)=\overline{X}\setminus X$ in $\overline{X}$, and let
$f_X:U_{b(X)}\to\R_+$ be a continuous function (constructed using
the control functions $\{\rho_X\}_{X\subseteq A}$ induced by the
compactification of $A$), smooth on the stratum $X$ such that
$f_X^{-1}(0)=b(X)$ and ${f_X}_{\vert U_{b(X)}\cap X}$ is
submersive. Now, apply lemma 4.1 to $v$; this gives a totally
radial vector field $v'$. Then we remark that if $v'$ is a totally
radial vector field, for all strata $X$, $v'$ is entering
$X_{\geq\epsilon}=X\setminus\{f_X<\epsilon\}$ along $\partial
X_{\epsilon}$ for $\epsilon$ small enough, where the symbol
$\partial X_{\epsilon}$ denotes the level hypersurface
$\{f_X=\epsilon\}$. This is because $grad(f_X)=\sum_{Y<X} a_Y.
grad(\rho_Y)$, where the $a_Y$ are non-negative smooth functions,
at every point of $X$. So we have
$\chi_c(X_{\geq\epsilon})-\chi_c(\partial X_{\epsilon})=
\sum_{v'(x)=0}Ind(v',x)$ thanks to the classical Poincar\'e-Hopf
theorem. Finally, we have $\chi^c(M)=\chi_c(M)-\chi_c(\partial M)$
for every compactifiable manifold $M$ by adding a boundary
$\partial M$. Use the ``additivity'' formula of theorem 3.1 and
the definition of the singular index to complete the proof.

\section{A few examples}
In the following computations, as we are only interested in the
rank of intersection homology groups, we shall take $R=\Q$ and
work with the dimension of $\Q-$vector spaces. Moreover, this will
permit us to apply Poincar\'e duality to calculate some associated
groups. In the remainder of the text, $T^2$ will denote the torus
$S^1\times S^1$. The stratifications of spaces will be the evident
ones and we shall not go into details. See \cite{Bo} for classical
tools to compute $IH_{\bullet}$ of the following spaces.
\subsection{An inevitable example : the pinched torus $T_p^2$}
We have a unique perversity
$\overline p=\overline 0$ and
we have evidently a totally radial vector field $u$ on $T_p^2$ with a unique
singularity at the isolated singular  point $x_0$ of $T_p^2$, of indice $1$.
The link at this point is $L_{x_0}=S^1\sqcup S^1$. We have~:
$$IH^{\overline 0}_i(T_p^2)=\cases{\Q & si $i=2$\cr
0 & si $i=1$\cr
\Q & si $i=0$}$$
\noindent
so that
$$I\chi^{\overline 0}(T_p^2)=2.$$
\noindent
On the other hand~:
$$\begin{array}{ll}
m^{\overline 0}_{x_0}(T_p^2)&=\dim IH_1^{\overline 0}(L_{x_0})\\
\null&=\dim H_1(S^1\sqcup S^1)\\
\null&=2.\\
\end{array}$$
\noindent
Finally we have
$I\chi^{\overline 0}(T_p^2)=2=2.1=
Ind^{\overline 0}(u,x_0).$
\subsection{A well-known example : the suspension of the torus $\Sigma T^2$
(H. Poincar\'e, {1895})}
This time, we have two different perversities $\overline 0$ and $\overline t$
and two isolated singularities (which are the two vertices of
suspension). The link at these points is $L_{x_0}=L_{x_1}=T^2$. We still have a
totally radial vector field $v$ with two singular points of indice $1$ at
singularities of $\Sigma T^2$.
Remark that this pseudo-manifold is normal so we have
$IH^{\overline t}_*(\Sigma T^2)=H_*(\Sigma T^2)$, i.e.
$$IH^{\overline t}_i(\Sigma T^2)=\cases{\Q & if $i=3$\cr
\Q^2 & if $i=2$\cr
0 & if $i=1$\cr
\Q & if $i=0$.}$$
\noindent
Hence
$$I\chi^{\overline t}(\Sigma T^2)=2$$
\noindent
and by duality we find
$$I\chi^{\overline 0}(\Sigma T^2)=-2.$$
\noindent
On the other hand~:
$$\begin{array}{ll}
m_{x_0}^{\overline 0}(\Sigma T^2)&=-\dim
IH^{\overline 0}_{2}(L_{x_0})\cr
\null&=-\dim H_2(T^2)\cr
\null&=-1
\end{array}$$
\noindent
and
$$\begin{array}{ll}m_{x_0}^{\overline t}(\Sigma T^2)&=\dim
IH^{\overline t}_{1}(L_{x_0})-\dim IH_2^{\overline t}(L_{x_0})\cr
\null&=\dim H_1(T^2)-\dim H_2(T^2)\cr
\null&=1.
\end{array}$$
\noindent
Finally we have~:
$$I\chi^{\overline 0}(\Sigma T^2)=-2=-1-1=
2.Ind^{\overline 0}(v,x_0)$$
\noindent
and
$$I\chi^{\overline t}(\Sigma T^2)=2=1+1=
2.Ind^{\overline t}(v,x_0).$$
\subsection{A hybrid example  : the suspension of the torus of dimension 3,
twice pinched, $\Sigma T^3_{2p}$}

We have
$\Sigma T^3_{2p}=\Sigma(\Sigma(T^2\sqcup T^2))$.
Here we have four perversities $\overline 0,\overline m,\overline n,
\overline t$.
\smallskip
\noindent
Calculate to begin with the homology of $T^3_{2p}$ :
$$H_i(T^3_{2p})=\cases{\Q^2 & si $i=3$\cr
\Q^4 & if $i=2$\cr
\Q & if $i=1$\cr
\Q & if $i=0$}.$$
\noindent
Then its intersection homology is :
$$\begin{array}{ll}
IH_i^{\overline 0}(T^3_{2p})&=\cases{%
H_i(T^3_{2p}) & if $i>2$\cr
Im\ (H_i(T^2\sqcup T^2)\rightarrow H_i(T^3_{2p})) & if $i=2$\cr
H_i(T^2\sqcup T^2) & if $i<2$}\cr
\null &=\cases{
\Q^2 & if $i=3$\cr
0 & if $i=2$\cr
\Q^4 & if $i=1$\cr
\Q^2 & if $i=0$}
\end{array}$$
\noindent
where we deduce
$$IH_i^{\overline t}(T^3_{2p})=\cases{%
\Q^2 & if $i=3$\cr
\Q^4 & if $i=2$\cr
0 & if $i=1$\cr
\Q^2 & if $i=0$}.$$
\noindent
And at last the intersection homology of the suspension $\Sigma T^3_{2p}$ is :
$$\begin{array}{ll}
IH^{\overline p}_i(\Sigma T^3_{2p})&=\cases{IH_{i-1}^{\overline p}(T^3_{2p}) & if $i>3-p_4$\cr
                                            0 & if $i=3-p_4$\cr
                                            IH_i^{\overline p}(T^3_{2p})  & if $i<3-p_4$}\cr

\null                              &=\cases{
\cases{
\Q^2 & if $i=4$\cr
0 & if $i=3$\cr
0 & if $i=2$\cr
\Q^4 & if $i=1$\cr
\Q^2 & if $i=0$} & if $\overline p=\overline m$\cr
\cases{
\Q^2 & if $i=4$\cr
0 & if $i=3$\cr
0 & if $i=2$\cr
\Q^4 & if $i=1$\cr
\Q^2 & if $i=0$} & if $\overline p=\overline 0.$%cases
}%cases
\end{array}$$
It is easy to construct a totally radial vector field $w$  with four
singularities~: two at the vertices of suspension, say $x_0,x_1$, of indice
$1$ and two others on
strata of codimension 3, say $x_2,x_3$, of indice $-1$.
Links are $L_{x_0}=L_{x_1}=T^3_{2p}$ and $L_{x_2}=L_{x_3}=T^2\sqcup T^2$.
Calculations of multiplicities give~:
$$m^{\overline p}_{x_0}(\Sigma T^3_{2p})=\cases{
\dim IH_1^{\overline t}(T^3_{2p})-\dim IH_2^{\overline t}(T^3_{2p})+
\dim IH_3^{\overline t}(T^3_{2p})=-2 & if $\overline p=\overline t$\cr
-\dim IH_2^{\overline t}(T^3_{2p})+\dim IH_3^{\overline t}(T^3_{2p})=-2 & if
$\overline p=\overline n$\cr
-\dim IH_2^{\overline 0}(T^3_{2p})+\dim IH_3^{\overline 0}(T^3_{2p})=2 & if
$\overline p=\overline m$\cr
\dim IH_3^{\overline 0}(T^3_{2p})=2 & if
$\overline p=\overline 0$}$$
\noindent
and
$$m^{\overline p}_{x_2}(\Sigma T^3_{2p})=\cases{
-\dim H_1(T^2\sqcup T^2)+\dim H_2(T^2\sqcup T^2)
=-2 & if $\overline p=\overline t$\cr
-\dim H_1(T^2\sqcup T^2)+\dim H_2(T^2\sqcup T^2)
=-2 & if $\overline p=\overline n$\cr
\dim H_2(T^2\sqcup T^2)=2 & if
$\overline p=\overline m$\cr\dim H_2(T^2\sqcup T^2)=2 & if
$\overline p=\overline 0$}.$$
\noindent
Finally we have
$$\begin{array}{ll}
I\chi^{\overline 0}(\Sigma T^3_{2p})&=0=2+2+2.(-1)+2.(-1)\cr
I\chi^{\overline m}(\Sigma T^3_{2p})&=0=2+2-2-2\cr
I\chi^{\overline n}(\Sigma T^3_{2p})&=0=-2-2+(-2).(-1)+(-2).(-1)\cr
I\chi^{\overline t}(\Sigma T^3_{2p})&=0=-2-2+2+2.
\end{array}$$
\goodbreak
\section{A partial converse}

We present here a partial converse to theorem 5.1 in the sense
that we study when a stratified set admits a strongly totally
radial vector field without singularity. This result is in the
line of \cite{Su}, \cite{V1} or \cite{Sc5}, \cite{Sc6}. See also
\cite{Mat2}, theorem 8.5. The result is partial because of the
example below. Indeed, it shows that we cannot expect the
condition $I\chi^{\overline p}_c(A)=0$ to imply the existence of a
non singular totally radial vector field.
\begin{theo} Let $A$ be a compactifiable $n$-pseudo-manifold.
There exists a strongly totally radial vector field
(relatively to $\cal A$) on $A$
without singularity if and only if $\chi^c(X)=0$ for all strata
$X$ of $\cal A$.
\end{theo}
{\it Proof}. To show sufficiency, we use the carpeting functions
of the proof of theorem 5.1. Let $v$ be a strongly totally radial
vector field on $A$ with isolated singularities~; the vector field
$v_X$ is entering on the boundary $\partial X_{\geq\epsilon}$
(defined by a level hypersurface of a carpeting function). Remark
that, as $\chi^c(X)=0$, we can deform $v_X$ on $X_{\geq\epsilon}$
(without modifying it near $\partial X_{\geq\epsilon}$) so as to
have no singularities {(\cite{H})}. We have evidently
${\rho_X}_*(v)>0$ on $T_X\setminus X$ for all strata $X$.
Necessity is proved in an analogous manner.
\begin{corol} Let $A$ be a compactifiable $n$-pseudo-manifold,
stratified with
strata of odd dimension. Then there exists a strongly totally radial vector
field
without singularity on $A$.
\end{corol}
{\it Remark}.
Existence of a totally radial vector
field
without singularity, on an abstract stratified set, is equivalent to the existence of a
controlled vector field without singularity.

\bigskip
\noindent
{\it Example.} Finally, here is an example of a compact pseudo-manifold without strata of
dimension
0 for which
$I\chi^{\overline p}_c(A)=0$ for every
%({\sl loose})
perversity $\overline p$ and admitting
no totally radial vector field without a singularity. Consider
$A=\Sigma (T_{2p}^3)\times S^2$~; it is clear that
$I\chi^{\overline p}(A;R)=0$ for all $\overline p$. Nevertheless, there does not
exist a totally radial vector field without a singularity (look at strata
$\{*\}\times S^2$ or $\{**\}\times S^2$). This is  also evident as a
consequence of theorem 7.1.
{\footnotesize
\bibliography{Poincare-HopfTypeHal}
}
\end{document}